\documentclass[12pt,a4paper]{article}
\usepackage{amssymb}

\usepackage{amsmath}

\input{tcilatex}

\begin{document}

\title{A Problem of Hsiang-Palais-Terng on Isoparametric Submanifolds}
\author{Haibao Duan \\
Institute of Mathematics, \\
Chinese Academy of Sciences,\\
Beijing 100080, dhb@math.ac.cn}
\date{ \ \ }
\maketitle

\begin{abstract}
We solve the problem raised by Hsiang, Palais and Terng in [HPT]: Is it
possible to have an isoparametic foliation on $\mathbb{R}^{52}$ whose marked
Dynkin diagram is of type $D_{4}$ and with all multiplicities uniformly
equal to $4$?

\begin{description}
\item \textsl{2000 Mathematical Subject classification: }53C42 (14M15; 57R20)

\item \textsl{Key words and phrases:} Isoparametric submanifolds; Vector
bundles; Characteristic classes.
\end{description}
\end{abstract}

\section{Introduction}

We assume familiarity with the notations, terminology and results developed
by Hsiang, Palais and Terng in [HPT].

Associated to an isoparametric submanifold $M$ in an Euclidean $n$-space $%
\mathbb{R}^{n}$, there is a Weyl group with Dynkin diagram marked with
multiplicities. In this paper we show

\textbf{Theorem.} \textsl{There is no isoparametric submanifold in }$\mathbb{%
R}^{52}$\textsl{\ whose marked Dynkin diagram is of type }$D_{4}$\textsl{\
and with all multiplicities uniformly equal to }$4$\textsl{.}

\bigskip

This result solves Problem 1 raised by Hsiang, Palais and Terng in [HPT]. As
was pointed out by the authors of [HPT], it implies also that

\textbf{Corollary.} \textsl{There is no isoparametric submanifolds whose
marked Dynkin diagrams with uniform multiplicity }$4$\textsl{\ of }$D_{k}$%
\textsl{-type, }$k>5$\textsl{\ or }$E_{k}$\textsl{-type, }$k=6,7,8$\textsl{.}

\section{Isoparametric submanifolds of type $D_{4}$}

Assume throughout this section that $M\subset \mathbb{R}^{n}$ is an
irreducible isoparametric submanifold with uniform even multiplicity $m$ of $%
D_{4}$-type. With this assumption we have

$\qquad \dim M=12m$ and $n=12m+4$.

\noindent The inner product on $\mathbb{R}^{n}$\ will be denoted by $(,)$.

Fix a base point $a\in M$ and let $P\subset \mathbb{R}^{n}$ be the normal
plane to $M$ at $a$. It is a subvector space with $\dim P=4$.

Let $\Lambda $ be the focal set of the embedding $M\subset \mathbb{R}^{n}$. $%
P$ intersects $\Lambda $ at $12$ linear hyperplanes (through the origin in $%
P $)

$\qquad \qquad \Lambda \cap P=L_{1}\cup \cdots \cup L_{12}$.

\noindent The reflection $\sigma _{i}$ of $P$ in the $L_{i}$, $1\leq i\leq
12 $, generate the Weyl group $W$ of type $D_{4}$, considered as a subgroup
of the isometries of $P$.

\bigskip

It follows that we can furnish $P$ with an orthonormal basis $\{\varepsilon
_{1},\varepsilon _{2},\varepsilon _{3},\varepsilon _{4}\}$ so that the set $%
\Phi =\{\alpha _{i}\in P\mid 1\leq i\leq 12\}$ of positive roots of $W$
relative to $a\in P$ is given and ordered by (cf. [Hu, p.64])

$\qquad \alpha _{1}=\varepsilon _{1}-\varepsilon _{2}$;$\qquad \alpha
_{2}=\varepsilon _{2}-\varepsilon _{3}$;$\qquad \alpha _{3}=\varepsilon
_{3}-\varepsilon _{4}$;

$\qquad \alpha _{4}=\varepsilon _{1}-\varepsilon _{3}$;$\qquad \alpha
_{5}=\varepsilon _{2}-\varepsilon _{4}$;$\qquad \alpha _{6}=\varepsilon
_{1}-\varepsilon _{4}$

$\qquad \alpha _{7}=\varepsilon _{2}+\varepsilon _{1}$;$\qquad \alpha
_{8}=\varepsilon _{3}+\varepsilon _{2}$;$\qquad \alpha _{9}=\varepsilon
_{4}+\varepsilon _{3}$;

$\qquad \alpha _{10}=\varepsilon _{3}+\varepsilon _{1}$;$\quad \ ~\alpha
_{11}=\varepsilon _{4}+\varepsilon _{2}$;$\quad \ \alpha _{12}=\varepsilon
_{4}+\varepsilon _{1}$.

\noindent Consequently, we order the planes $L_{i}$ by the requirement that
each $\alpha _{i}$ is normal to $L_{i}$, $1\leq i\leq 12$. As results we have

\noindent (2-1) $\alpha _{1}$, $\alpha _{2},\alpha _{3}$ and $\alpha _{9}$
form the set of simple roots relative to $a$, and the corresponding Cartan
matrix is

$\qquad (\beta _{ij})_{i,j=1,2,3,9}=\left(
\begin{array}{cccc}
2 & -1 & 0 & 0 \\
-1 & 2 & -1 & -1 \\
0 & -1 & 2 & 0 \\
0 & -1 & 0 & 2%
\end{array}%
\right) $; $\beta _{ij}=2\frac{(\alpha _{i},\alpha _{j})}{(\alpha
_{j},\alpha _{j})}$.

\noindent (2-2) The group $W$ is generated $\sigma _{1}$, $\sigma
_{2},\sigma _{3}$ and $\sigma _{9}$, whose actions on $P$ are given
respectively by

$\qquad \sigma _{1}:\{\varepsilon _{1},\varepsilon _{2},\varepsilon
_{3},\varepsilon _{4}\}\rightarrow \{\varepsilon _{2},\varepsilon
_{1},\varepsilon _{3},\varepsilon _{4}\}$,

$\qquad \sigma _{2}:\{\varepsilon _{1},\varepsilon _{2},\varepsilon
_{3},\varepsilon _{4}\}\rightarrow \{\varepsilon _{1},\varepsilon
_{3},\varepsilon _{2},\varepsilon _{4}\}$,

$\qquad \sigma _{3}:\{\varepsilon _{1},\varepsilon _{2},\varepsilon
_{3},\varepsilon _{4}\}\rightarrow \{\varepsilon _{1},\varepsilon
_{2},\varepsilon _{4},\varepsilon _{3}\}$,

$\qquad \sigma _{9}:\{\varepsilon _{1},\varepsilon _{2},\varepsilon
_{3},\varepsilon _{4}\}\rightarrow \{\varepsilon _{1},\varepsilon
_{2},-\varepsilon _{4},-\varepsilon _{3}\}$;

\noindent (2-3) Let $b=\varepsilon _{1}+\varepsilon _{2}+\varepsilon
_{3}+\varepsilon _{4}\in P$. The subgroup $W_{b}$ of $W$ that fixes $b$ is
generated by $\sigma _{1}$, $\sigma _{2}$ and $\sigma _{3}$;

Using (2-2) one verifies directly that

\noindent (2-4) the set simple roots can be expressed in term of the $W$%
-action on $\alpha _{1}$ as

$\qquad \alpha _{2}=\sigma _{1}\sigma _{2}(\alpha _{1})$; $\qquad \alpha
_{3}=\sigma _{2}\sigma _{1}\sigma _{3}\sigma _{2}(\alpha _{1})$; $\ \alpha
_{4}=\sigma _{2}(\alpha _{1})$;

$\qquad \alpha _{5}=\sigma _{1}\sigma _{3}\sigma _{2}(\alpha _{1})$; $\quad
\alpha _{6}=\sigma _{3}\sigma _{2}(\alpha _{1})$;

$\qquad \alpha _{7}=\sigma _{2}\sigma _{3}\sigma _{9}\sigma
_{2}(\alpha _{1}) $; $\alpha _{8}=\sigma _{1}\sigma _{3}\sigma
_{9}\sigma _{2}(\alpha _{1})$; $\quad \alpha _{9}=\sigma
_{2}\sigma _{1}\sigma _{9}\sigma _{2}(\alpha _{1}) $;

$\qquad \alpha _{10}=\sigma _{3}\sigma _{9}\sigma _{2}(\alpha _{1})$; $\ \
\alpha _{11}=\sigma _{1}\sigma _{9}\sigma _{2}(\alpha _{1})$; $\quad \alpha
_{12}=\sigma _{9}\sigma _{2}(\alpha _{1})$.

\textbf{Remark 1.} (2-1) implies the following geometric facts. The planes $%
L_{i}$ partition $P$ into $\mid W\mid =2^{3}\cdot 4!$ convex open hulls
(called \textsl{Weyl chambers}), and the base point $a$ is contained in the
one $\Omega $ bounded by the $L_{i}$, $1\leq i\leq 3$ and $L_{9}$. In (2-3)
the point $b$ lies on the edge $L_{1}\cap L_{2}\cap L_{3}$ of $\overline{%
\Omega }$.$\square $

\bigskip

let $M_{b}\subset \mathbb{R}^{n}$ be the focal manifold parallel to $M$
through $b$ [HPT]. We have a smooth projection

$\qquad \qquad \pi :M\rightarrow M_{b}$

\noindent whose fiber over $b\in M_{b}$ is denoted by $F$.

Recall from [HPT] that the tangent bundle $TM$ of $M$ has a canonical
splitting as the orthogonal direct sum of $12$ subbundles

$\qquad \qquad TM=E_{\alpha _{1}}\oplus \cdots \oplus E_{\alpha _{12}}$, $%
\dim _{\mathbb{R}}E_{\alpha _{i}}=m$,

\noindent in which $E_{i}$ is \textsl{the curvature distribution} of $M$
relative to the root $\alpha _{i}$, $1\leq i\leq 12$. From [HPT] we have

\textbf{Lemma 1.} \textsl{Let }$TN$\textsl{\ be the tangent bundle of a
smooth manifold} $N$\textsl{. Then}

\textsl{(1) the subbundle }$\oplus _{1\leq i\leq 6}E_{\alpha _{i}}$\textsl{\
of }$TM$\textsl{\ restricts to }$TF$\textsl{; }

\textsl{(2) the induced bundle }$\pi ^{\ast }TM_{b}$\textsl{\ agrees with }$%
\oplus _{7\leq i\leq 12}E_{\alpha _{i}}$.$\square $

\section{The cohomology of the fibration $\protect\pi :$ $M\rightarrow M_{b}$%
}

Let $b_{i}\in H_{m}(M;\mathbb{Z})$, $1\leq i\leq 12$, be the homology class
of the leaf sphere $S_{i}(a)\subset M$ of the intergrable bundle $E_{\alpha
_{i}}$ through $a\in M$, and let $d_{i}\in H^{m}(M;\mathbb{Z})$ be the Euler
class of $E_{\alpha _{i}}$.

\textbf{Lemma 2.} \textsl{The Kronecker pairing }$<,>:H^{m}(M;\mathbb{Z}%
)\otimes H_{m}(M;\mathbb{Z})\rightarrow \mathbb{Z}$ \textsl{can be expressed
in term of the inner product }$(,)$\textsl{\ on }$P$\textsl{\ as}

$\qquad \qquad <d_{i},b_{j}>=2\frac{(\alpha _{i},\alpha _{j})}{(\alpha
_{j},\alpha _{j})}$\textsl{, }$1\leq i,j\leq 12$\textsl{.}$\square $

\textbf{Remark 2.} Let $G$ be a compact connected semi-simple Lie group with
a fixed maximal torus $T$ and Weyl group $W$. Fix a regular point $a$ in the
Cartan subalgebra $L(T)$ of the Lie algebra $L(G)$ corresponding to $T$. The
orbit of the adjoint action of $G$ on $L(G)$ through $a$ yields an embedding
$G/T\rightarrow L(G)$ which defines the flag manifold $G/T$ as an
isoparametric submanifold in $L(G)$ with associated Weyl group $W$ and with
equal multiplicities $m=2$.

In this case Lemma 2 has its generality due to Bott and Samelson [BS], and
the numbers $2\frac{(\alpha _{i},\alpha _{j})}{(\alpha _{j},\alpha _{j})}$
are \textsl{the Cartan numbers} of $G$ (only $0$, $\pm 1,\pm 2,\pm 3$ can
occur).$\square $

\bigskip

Since the roots $\alpha _{1}$,$\alpha _{2},\alpha _{3}$ and $\alpha _{9}$
form a set of simple roots, the classes $b_{1}$,$b_{2},b_{3}$ and $b_{9}$
constitute an additive basis of $H_{m}(M;\mathbb{Z})$. Since the $M$ is $m-1$
connected, we specify a basis $\omega _{1},\omega _{2},\omega _{3},\omega
_{9}$ of $H^{m}(M;\mathbb{Z})$ (in term of Kronecker pairing) as

$\qquad <b_{i},\omega _{j}>=\delta _{ij}$, $i,j=1,2,3,9$.

\noindent It follows from Lemma 2 that

\textbf{Lemma 3.} $\left(
\begin{array}{c}
d_{1} \\
d_{2} \\
d_{3} \\
d_{9}%
\end{array}%
\right) =\left(
\begin{array}{cccc}
2 & -1 & 0 & 0 \\
-1 & 2 & -1 & -1 \\
0 & -1 & 2 & 0 \\
0 & -1 & 0 & 2%
\end{array}%
\right) \left(
\begin{array}{c}
\omega _{1} \\
\omega _{2} \\
\omega _{3} \\
\omega _{9}%
\end{array}%
\right) .\square $

\bigskip

In term of the $\omega _{i}$ we introduce in $H^{m}(M;\mathbb{Z})$ the
classes $t_{i}$ , $1\leq i\leq 4$, by the relation

\noindent (3-1)$\qquad \qquad \left(
\begin{array}{c}
t_{1} \\
t_{2} \\
t_{3} \\
t_{4}%
\end{array}%
\right) =\left(
\begin{array}{cccc}
1 & 0 & 0 & 0 \\
-1 & 1 & 0 & 0 \\
0 & -1 & 1 & 0 \\
0 & 0 & -1 & 2%
\end{array}%
\right) \left(
\begin{array}{c}
\omega _{1} \\
\omega _{2} \\
\omega _{3} \\
\omega _{9}%
\end{array}%
\right) $.

\noindent Conversely,

\noindent (3-2)$\qquad \qquad \left(
\begin{array}{c}
\omega _{1} \\
\omega _{2} \\
\omega _{3} \\
\omega _{9}%
\end{array}%
\right) =\left(
\begin{array}{cccc}
1 & 0 & 0 & 0 \\
1 & 1 & 0 & 0 \\
1 & 1 & 1 & 0 \\
\frac{1}{2} & \frac{1}{2} & \frac{1}{2} & \frac{1}{2}%
\end{array}%
\right) \left(
\begin{array}{c}
t_{1} \\
t_{2} \\
t_{3} \\
t_{4}%
\end{array}%
\right) $.

\bigskip

The Weyl group $W$ (acting as isometries of $P$) has the effect to permute
roots (cf. (2-3)). On the other hand, $W$ acts also smoothly on $M$ [HPT],
hence acts as automorphisms of the cohomology (resp. the homology) of $M$.
For an $w\in W$ write $w^{\ast }$ (resp. $w_{\ast }$) for the induced action
on the cohomology (resp. homology).

\textbf{Lemma 4.}\textsl{\ With respect to the }$\mathbb{Q}$\textsl{\ basis }%
$t_{1},t_{2},t_{3},t_{4}$\textsl{\ of }$H^{2}(M;\mathbb{Q})$\textsl{, the
action of }$W$\textsl{\ on }$H^{2}(M;\mathbb{Q})$\textsl{\ is given by}

$\qquad \sigma _{1}^{\ast }:\{t_{1},t_{2},t_{3},t_{4}\}\rightarrow
\{t_{2},t_{1},t_{3},t_{4}\},$

$\qquad \sigma _{2}^{\ast }:\{t_{1},t_{2},t_{3},t_{4}\}\rightarrow
\{t_{1},t_{3},t_{2},t_{4}\},$

$\qquad \sigma _{3}^{\ast }:\{t_{1},t_{2},t_{3},t_{4}\}\rightarrow
\{t_{1},t_{2},t_{4},t_{3}\}$

$\qquad \sigma _{9}^{\ast }:\{t_{1},t_{2},t_{3},t_{4}\}\rightarrow
\{t_{1},t_{2},-t_{4},-t_{3}\}$.

\textbf{Proof. }Let $i,j,k=1,2,3,9$. In term of the Cartan matrix (2-1), the
action of $\sigma _{i\ast }$ on the $\mathbb{Z}$-basis $%
b_{1},b_{2},b_{3},b_{9}$ has been determined in [HPT] as

\noindent (3-3)$\qquad \sigma _{i\ast }(b_{j})=b_{j}+\beta _{ij}b_{i}$

\noindent (cf. the proof of 6.11. Corollary in [HPT]). By the naturality of
Kronecker pairing

$\qquad <\sigma _{i}^{\ast }(\omega _{k}),b_{i}>=<\omega _{k},\sigma _{i\ast
}(b_{j})>=\delta _{kj}+\beta _{ij}\delta _{ki}$,

\noindent we get from (3-3) that

\noindent (3-4) $\quad \sigma _{i}^{\ast }(\omega _{k})=\{%
\begin{array}{c}
\omega _{k}\text{ if }i\neq k\text{;\qquad \qquad \qquad \qquad \qquad
\qquad \qquad \qquad } \\
\omega _{k}-(\beta _{k1}\omega _{1}+\beta _{k2}\omega _{2}+\beta _{k3}\omega
_{3}+\beta _{k9}\omega _{9})\text{ if }i=k\text{.}%
\end{array}%
$

\noindent With $(\beta _{ij})$ being given explicitly in (2-1),
combining (3-1), (3-4) with (3-2) verifies Lemma 4.$\square $

\bigskip

Let the algebra $\mathbb{Q}[t_{1},t_{2},t_{3},t_{4}]$ of polynomials in the
variables $t_{1},t_{2},t_{3},t_{4}$ be graded by $\deg (t_{i})=m$, $1\leq
i\leq 4$. Let $e_{i}\in \mathbb{Q}[t_{1},t_{2},t_{3},t_{4}]$ (resp. $\theta
_{i}\in \mathbb{Q}[t_{1},t_{2},t_{3},t_{4}]$) be the $i^{th}$ elementary
symmetric functions in $t_{1},t_{2},t_{3},t_{4}$ (resp. in $%
t_{1}^{2},t_{2}^{2},t_{3}^{2},t_{4}^{2}$), $1\leq i\leq 4$. We note that the
$\theta _{i}$ can be written as a polynomial in the $e_{i}$

$\qquad \theta _{1}=e_{1}^{2}-2e_{2}$;$\quad \theta
_{2}=e_{2}^{2}-2e_{1}e_{3}+2e_{4}$;$\quad \theta _{3}=e_{3}^{2}-2e_{2}e_{4}$.

\textbf{Lemma 5. }\textsl{In term of generator-relations the rational
cohomology of }$M$\textsl{\ is given by}

$\qquad H^{\ast }(M;\mathbb{Q})=\mathbb{Q}[t_{1},t_{2},t_{3},t_{4}]/\mathbb{Q%
}^{+}[\theta _{1},\theta _{2},\theta _{3},e_{4}]$\textsl{.}

\noindent \textsl{Further, the induced homomorphism }$\pi ^{\ast }:H^{\ast
}(M_{b};\mathbb{Q})\rightarrow H^{\ast }(M;\mathbb{Q})$\textsl{\ maps the
algebra }$H^{\ast }(M_{b};\mathbb{Q})$\textsl{\ isomorphically onto the
subalgebra }

$\qquad \qquad\mathbb{Q}[e_{1},e_{2},e_{3}]/\mathbb{Q}^{+}[\theta
_{1},\theta _{2},\theta _{3}]$\textsl{.}

\textbf{Proof.} It were essentially shown in [HPT, 6.12. Theorem; 6.14.
Theorem] that for any $\mathbb{Q}$ basis $y_{1},y_{2},y_{3},y_{4}$ of $%
H^{m}(M;\mathbb{Q})$ one has the grade preserving $W$-isomorphisms

$\qquad H^{\ast }(M;\mathbb{Q})=\mathbb{Q}[y_{1},y_{2},y_{3},y_{4}]/\mathbb{Q%
}^{+}[y_{1},y_{2},y_{3},y_{4}]^{W}$;

$\qquad H^{\ast }(M_{b};\mathbb{Q})=\mathbb{Q}%
[y_{1},y_{2},y_{3},y_{4}]^{W_{b}}/\mathbb{Q}%
^{+}[y_{1},y_{2},y_{3},y_{4}]^{W} $,

\noindent where the $W$-action on the $\mathbb{Q}$-algebra $\mathbb{Q}%
[y_{1},y_{2},y_{3},y_{4}]$ is induced from the $W$-action on the $\mathbb{Q}$%
-vector space $H^{m}(M;\mathbb{Q})=span_{Q}\{y_{1},y_{2},y_{3},y_{4}\}$ and
where $\mathbb{Q}[y_{1},y_{2},y_{3},y_{4}]^{W_{b}}$ (resp. $\mathbb{Q}%
^{+}[y_{1},y_{2},y_{3},y_{4}]^{W}$) is the subalgebra of $W_{b}$-invariant
polynomials (resp. $W$-invariant polynomials in positive degrees).

Since the transition matrix from $t_{1},t_{2},t_{3},t_{4}$ to the $\mathbb{Z}
$-basis $\omega _{1},\omega _{2},\omega _{3},\omega _{9}$ of $H^{m}(M;%
\mathbb{Z})$ is non-singular by (3-1), the $t_{i}$ constitute a basis for $%
H^{m}(M;\mathbb{Q})$. Moreover, it follows from Lemma 4 that

$\mathbb{Q}[t_{1},t_{2},t_{3},t_{4}]^{W_{b}}=\mathbb{Q}%
[e_{1},e_{2},e_{3},e_{4}]$,$\quad \mathbb{Q}[t_{1},t_{2},t_{3},t_{4}]^{W}=%
\mathbb{Q}[\theta _{1},\theta _{2},\theta _{3},e_{4}]$.

\noindent This completes the proof.$\square $

\textbf{Remark 3.} If $G=SO(2n)$ (the special orthogonal group of rank $2n$%
), the embedding $G/T\rightarrow L(G)$ considered in Remark 2 gives rise to
an isoparametric submanifold $M=SO(2n)/T$ in the Lie algebra $L(SO(2n))$
which is of $D_{n}$-type with equal multiplicities $m=2$. The corresponding $%
M_{b}$ is known as the \textsl{Grassmannian of complex structures} on $%
\mathbb{R}^{2n}$ [D]. In this case Borel computed the algebras $H^{\ast }(M;%
\mathbb{Q})$ and $H^{\ast }(M_{b};\mathbb{Q})$ in [B] which are compatible
with Lemma 5.$\square $

\section{Computation in the Pontrijagin classes}

Turn to the case $m=4$ concerned by our Theorem. Denote by $p_{1}(\xi )\in
H^{4}(X;\mathbf{Z})$ for the first Pontrijagin class of a real vector bundle
$\xi $ over a topological space $X$.

\textbf{Lemma 6.} \textsl{For an }$w\in W$\textsl{\ and an }$\alpha \in \Phi
$ \textsl{one has}

\textsl{\qquad \qquad }$w^{\ast }(p_{1}(E_{\alpha }))=p_{1}(E_{w^{-1}(\alpha
)})$\textsl{.}

\noindent \textsl{In particular, (2-3) implies that}

$p_{1}(E_{\alpha _{2}})=\sigma _{1}^{\ast }\sigma _{2}^{\ast
}p_{1}(E_{\alpha _{1}})$; $\qquad p_{1}(E_{\alpha _{3}})=\sigma _{2}^{\ast
}\sigma _{1}^{\ast }\sigma _{3}^{\ast }\sigma _{2}^{\ast }p_{1}(E_{\alpha
_{1}})$;

$p_{1}(E_{\alpha _{4}})=\sigma _{2}^{\ast }p_{1}(E_{\alpha _{1}})$; $\quad
\quad \quad p_{1}(E_{\alpha _{5}})=\sigma _{1}^{\ast }\sigma _{3}^{\ast
}\sigma _{2}^{\ast }p_{1}(E_{\alpha _{1}})$;

$p_{1}(E_{\alpha _{6}})=\sigma _{3}^{\ast }\sigma _{2}^{\ast
}p_{1}(E_{\alpha _{1}})$; $\quad \quad p_{1}(E_{\alpha _{7}})=\sigma
_{2}^{\ast }\sigma _{3}^{\ast }\sigma _{9}^{\ast }\sigma _{2}^{\ast
}p_{1}(E_{\alpha _{1}})$;

$p_{1}(E_{\alpha _{8}})=\sigma _{1}^{\ast }\sigma _{3}^{\ast }\sigma
_{9}^{\ast }\sigma _{2}^{\ast }p_{1}(E_{\alpha _{1}})$; $p_{1}(E_{\alpha
_{9}})=\sigma _{2}^{\ast }\sigma _{1}^{\ast }\sigma _{9}^{\ast }\sigma
_{2}^{\ast }p_{1}(E_{\alpha _{1}})$;

$p_{1}(E_{\alpha _{10}})=\sigma _{3}^{\ast }\sigma _{9}^{\ast }\sigma
_{2}^{\ast }p_{1}(E_{\alpha _{1}}))$; $\ p_{1}(E_{\alpha _{11}})=\sigma
_{1}^{\ast }\sigma _{9}^{\ast }\sigma _{2}^{\ast }p_{1}(E_{\alpha _{1}})$;

$p_{1}(E_{\alpha _{12}})=\sigma _{9}^{\ast }\sigma _{2}^{\ast
}p_{1}(E_{\alpha _{1}})$.

\textbf{Proof.} In term of the $W$-action on the set of roots, the induced
bundle $w^{\ast }(E_{\alpha _{i}})$ is $E_{w^{-1}(\alpha _{i})}$ (cf. [HTP,
\textbf{1.6}]). Lemma 6 comes now from the naturality of Pontrijagin classes
and from (2-4).$\square $

\bigskip

\textbf{Lemma 7.} $p_{1}(E_{\alpha _{1}})=k(t_{1}+t_{2}-t_{3}-t_{4})$
\textsl{for some} $k\in \mathbb{Q}$.

\textbf{Proof.} In view of Lemma 4 we can assume that

$\qquad \qquad p_{1}(E_{\alpha _{1}})=k_{1}t_{1}+\cdots +k_{4}t_{4}$, $%
k_{i}\in \mathbb{Q}$.

\noindent Since the restricted bundle $E_{\alpha _{1}}\mid S_{1}(a)$ is the
tangent bundle of the $4$-sphere $S_{1}(a)$ and therefore is stably trivial,
we have

$\qquad <p_{1}(E_{\alpha _{1}}),b_{1}>=k_{1}-k_{2}=0$.

\noindent That is

\noindent (4-1) $\qquad p_{1}(E_{\alpha
_{1}})=kt_{1}+kt_{2}+k_{3}t_{3}+k_{4}t_{4}$.

\noindent Since the actions of the $\sigma _{i}^{\ast }$, $i=1,2,3,9,$ on
the $t_{j}$ are known by Lemma 4, combining (4-1) with the relations in
Lemma 6 yields

\noindent (4-2)\qquad \qquad $%
\begin{array}{c}
p_{1}(E_{\alpha _{2}})=k_{3}t_{1}+kt_{2}+kt_{3}+k_{4}t_{4}; \\
p_{1}(E_{\alpha _{3}})=k_{3}t_{1}+k_{4}t_{2}+kt_{3}+kt_{4}; \\
p_{1}(E_{\alpha _{4}})=kt_{1}+k_{3}t_{2}+kt_{3}+k_{4}t_{4}; \\
p_{1}(E_{\alpha _{5}})=k_{3}t_{1}+kt_{2}+k_{4}t_{3}+kt_{4}; \\
p_{1}(E_{\alpha _{6}})=kt_{1}+k_{3}t_{2}+k_{4}t_{3}+kt_{4}; \\
p_{1}(E_{\alpha _{7}})=kt_{1}-kt_{2}+k_{3}t_{3}-k_{4}t_{4}; \\
p_{1}(E_{\alpha _{8}})=k_{3}t_{1}+kt_{2}-kt_{3}-k_{4}t_{4}; \\
p_{1}(E_{\alpha _{9}})=k_{3}t_{1}-k_{4}t_{2}+kt_{3}-kt_{4}; \\
p_{1}(E_{\alpha _{10}})=kt_{1}+k_{3}t_{2}-kt_{3}-k_{4}t_{4}; \\
p_{1}(E_{\alpha _{11}})=k_{3}t_{1}+kt_{2}-k_{4}t_{3}-kt_{4}; \\
p_{1}(E_{\alpha _{12}})=kt_{1}+k_{3}t_{2}-k_{4}t_{3}-kt_{4}\text{.}%
\end{array}%
$

Since the tangent bundle of any isoparametric submanifold is stably trivial,
we get from $TM=E_{\alpha _{1}}\oplus \cdots \oplus E_{\alpha _{12}}$ that

$\qquad p_{1}(TM)=p_{1}(E_{\alpha _{1}})+\cdots +p_{1}(E_{\alpha _{12}})=0$.

\noindent Comparing the coefficients of $t_{1}$ on both sides of the
equation turns out $k_{3}=-k$. Substituting this in (4-2) gives rise to, in
particular, that

\noindent (4-3)$\qquad
\begin{array}{c}
p_{1}(E_{\alpha _{7}})=kt_{1}-kt_{2}-kt_{3}-k_{4}t_{4};\ ~ \\
p_{1}(E_{\alpha _{8}})=-kt_{1}+kt_{2}-kt_{3}-k_{4}t_{4}; \\
p_{1}(E_{\alpha _{9}})=-kt_{1}-k_{4}t_{2}+kt_{3}-kt_{4}; \\
p_{1}(E_{\alpha _{10}})=kt_{1}-kt_{2}-kt_{3}-k_{4}t_{4};~ \\
p_{1}(E_{\alpha _{11}})=-kt_{1}+kt_{2}-k_{4}t_{3}-kt_{4}; \\
p_{1}(E_{\alpha _{12}})=kt_{1}-kt_{2}-k_{4}t_{3}-kt_{4}\text{.~~}%
\end{array}%
$

Finally, since

$ p_{1}(E_{\alpha _{7}})+\cdots +p_{1}(E_{\alpha _{12}})=\pi
^{\ast }p_{1}(TM_{b})\in$\textsl{\ Im } $[\pi ^{\ast }:
H^{\ast}(M_{b};\mathbb{Q})\rightarrow H^{\ast }(M;\mathbb{Q})]$

\noindent by (2) of Lemma 1, it must be symmetric in $t_{1},t_{2},t_{3},t_{4}
$ by (2) of Lemma 4. Consequently, $k_{4}=-k$. This completes the proof of
Lemma 6.$\square $

\bigskip

We emphasis what we actually need in the next result.

\textbf{Lemma 8.} \textsl{If }$M\subset R^{52}$\textsl{\ is an irreducible
isoparametric submanifold with uniform multiplicity }$4$\textsl{\ of }$D_{4}$%
\textsl{-type, there exists a }$4$\textsl{-plane bundle }$\xi $\textsl{\
over }$M$\textsl{\ whose Euler and the first Pontrijagin classes are
respectively}

\noindent (4-4)$\qquad \qquad e(\xi )=2\omega _{1}-\omega _{2}$;\ \quad $%
p_{1}(\xi )=2k(\omega _{2}-\omega _{4})$

\noindent \textsl{for some }$k\in \mathbb{Z}$\textsl{.}

\textbf{Proof.} Take $\xi =E_{\alpha _{1}}$. Then $e(\xi )=d_{1}=2\omega
_{1}-\omega _{2}$ by Lemma 3 and

$\qquad p_{1}(\xi )=k(t_{1}+t_{2}-t_{3}-t_{4})$ (by Lemma 7)

$\qquad \qquad \qquad =2k(\omega _{2}-\omega _{9})$ (by (3-1))

\noindent for some $k\in \mathbb{Q}$. Moreover we must have $k\in \mathbb{Z}$
since

(1) $M$ is $3$ connected and the classes $\omega _{1},\omega _{2},\omega
_{3},\omega _{9}$ constitute an additive basis \ for $H^{4}(M;\mathbb{Z})$;
and since

(2) the first Pontrijagin class of any vector bundle over a $3$ connected
CW-complex is an integer class and is divisible by $2$ [LD].$\square $

\section{A topological constraint on isoparametric submanifolds with equal
multiplicity $4$}

Let $Vect^{m}(S^{n})$ be the set of isomorphism classes of Euclidean $m$%
-vector bundles over the $n$-sphere $S^{n}=\{(x_{1},\cdots ,x_{n+1})\in
\mathbb{R}^{n+1}\mid x_{1}^{2}+\cdots +x_{n+1}^{2}=1\}$.

If $n=m=4$ we introduce the map $f:Vect^{4}(S^{4})\rightarrow \mathbb{Z}%
\oplus \mathbb{Z}$ by

$\qquad \qquad f(\xi )=(<e(\xi ),[S^{4}]>,<p_{1}(\xi ),[S^{4}]>),$

\noindent where

(1) $[S^{4}]\in H_{4}(S^{4};\mathbb{Z})=\mathbb{Z}$ is a fixed orientation
class;

(2) $<,>$ is the Kronecker pairing between cohomology and homology;

\noindent and where

(3) $e(\xi )$ and $p_{1}(\xi )$ are respectively the Euler and the first
Pontrijagin classes of $\xi \in Vect^{4}(S^{4})$.

\textbf{Example} (cf. [MS, p.246])\textbf{.} Let $\tau \in Vect^{4}(S^{4})$
be the tangent bundle of $S^{4}$, and let $\gamma \in Vect^{4}(S^{4})$ be
the real reduction of the quaternionic line bundle over $HP^{1}=S^{4}$ ($1$%
-dimensional quaternionic projective space). Then

$\qquad f(\tau )=(2,0);$ $\qquad f(\gamma )=(1,-2)$.

\bigskip

Our theorem will follow directly from Lemma 8 and the next result that
improves Lemma 20.10 in [MS].

\textbf{Lemma 9.} $f$\textsl{\ fits in the short exact sequence}

$\qquad \qquad 0\rightarrow Vect^{4}(S^{4})\overset{f}{\rightarrow
}\mathbb{Z}\oplus \mathbb{Z}\overset{g}{\mathbb{\rightarrow
}}\mathbb{Z} _{4}\rightarrow 0,$

\noindent \textsl{where} $g(a,b)\equiv (2a-b)$ \textsl{\ mod }$4$
\textsl{.}

\bigskip

\textbf{Proof of the Theorem.} Assume that there exists an irreducible
isoparametric submanifold $M\subset R^{52}$ with uniform multiplicity $4$ of
$D_{4}$-type. Let $\xi $ be a $4$-plane bundle over $M$ whose Euler and the
first Pontrijagin classes are given as that in (4-4).

Let $\xi _{1},\xi _{2}\in Vect^{4}(S^{4})$ be obtained respectively by
restricting $\xi $ to the leaf spheres $S_{2}(a)$, $S_{9}(a)$ (cf. section
3). By the naturality of characteristic classes we get from (4-4) that

$\qquad \qquad f(\xi _{1})=(-1,2k)$; $\qquad f(\xi _{2})=(0,-2k)$,

\noindent where the integer $k$ must satisfy the congruences

$\qquad \qquad -2\equiv 2k$ \textsl{\ mod }$4$; $\quad 0\equiv
-2k$ \textsl{\ mod }$4$.

\noindent by Lemma 9. The proof is done by the obvious
contradiction.$\square$

\bigskip

It suffices now to justify Lemma 9.

Let $SO(m)$ be the special orthogonal group of rank $m$ and denote by $\pi
_{r}(X)$ the $r$-homotopy group of \ a topological space $X$. The \textsl{%
Steenrod correspondence} is the map $s:Vect^{m}(S^{n})\rightarrow \pi
_{n-1}(SO(m))$ defined by

$s(\xi )=$the homotopy class of a clutching function $S^{n-1}\rightarrow
SO(m)$ of $\xi $.

\noindent In [S, \S 18], Steenrod showed that

\textbf{Lemma 10.} $s$\textsl{\ is a one-to-one correspondence.}$\square $

It follows that $Vect^{m}(S^{n})$ has a group structure so that $s$ is a
group isomorphism. We clarify this structure in Lemma 11.

\bigskip

Fix a base point $s_{0}=(1,0,\cdots ,0)\in S^{n}$. For two $\xi ,\eta \in
Vect^{m}(S^{n})$ write $\xi \vee \eta $ for the $m$-bundle over $S^{n}\vee
S^{n}$ (one point union of two $S^{n}$ over $s_{0}$) whose restriction to
the first (resp. the second) sphere agrees with $\xi $ (resp. $\eta $).
Define the addition $+:Vect^{m}(S^{n})\times Vect^{m}(S^{n})\rightarrow
Vect^{m}(S^{n})$ and the inverse $-:Vect^{m}(S^{n})\rightarrow
Vect^{m}(S^{n})$ operations by the rules

$\qquad \xi +\eta =\mu ^{\ast }(\xi \vee \eta )$ and $-\xi =\nu ^{\ast }\xi $%
,

\noindent where $\xi ,\eta \in Vect^{m}(S^{n})$, $\mu :S^{n}\rightarrow
S^{n}\vee S^{n}$ is the map that pinches the equator $x\in S^{n}$, $%
x_{n+1}=0 $ to the base point $s_{0}$ and where $\nu :S^{n}\rightarrow S^{n}$
is the restriction of the reflection of $\mathbb{R}^{n+1}$ in the hyperplane
$x_{n+1}=0$. It is straight forward to see that

\textbf{Lemma 11.} \textsl{With respect to the operations
}$+$\textsl{ and } $-$

\textsl{(1) }$Vect^{m}(S^{n})$\textsl{\ is an abelian group with
zero }$ \varepsilon ^{m}$\textsl{, the trivial }$m$\textsl{-bundle
over }$ S^{n}$\textsl{\ and}

\textsl{(2) the maps }$s$\textsl{\ and }$f$\textsl{\ are
homomorphisms.}$ \square $

\bigskip

We are ready to show Lemma 9.

\textbf{The proof of Lemma 9.} It is essentially shown by Milnor
[MS, p.245] that $f$ is injective and satisfies Im $f \supseteq$
Ker $g$. It suffices to show that $\textsl{\ Im }f\subseteq
\textsl{\ Ker }g$.

Assume on the contrary that there is a $\xi \in Vect^{4}(S^{4})$ such that $%
f(\xi )=(a,b)$ with $2a-b=4k+i$, $i=1,2,3$. Moreover, one must has $i=2$
since the first Pontrijagin class of any vector bundle over $S^{4}$ is
divisible by $2$. That is

\noindent (5-1) $\qquad f(\xi )=(a,2a-4k-2)$

Using the group operations in $Vect^{4}(S^{4})$ we form the class

$\qquad \widehat{\xi }=-\xi -(a-2k-1)\gamma +(a-k)\tau $,

\noindent where $\tau $ and $\gamma $ were given in the Example. By the
additivity of $f$ we get from the Example and (5-1) that

\noindent (5-2) $\qquad f(\widehat{\xi })=(1,0)$.

Consider the following diagram

\noindent $%
\begin{array}{ccc}
0\rightarrow \pi _{4}(S^{4})\overset{\partial }{\rightarrow } & \pi
_{3}(SO(4))=Vect^{4}(S^{4})\overset{i_{\ast }}{\rightarrow }\pi
_{3}(SO(5))=Vect^{5}(S^{4}) & \rightarrow 0 \\
& p_{1}\searrow \qquad \swarrow p_{1} &  \\
& H^{4}(S^{4};\mathbb{Z}) &
\end{array}%
$

\noindent in which

(1) the top row is a section in the homotopy exact sequence of the fibration
$SO(4)\overset{i}{\hookrightarrow }SO(5)\rightarrow S^{4}$ (cf. [Wh, p.196]);

(2) $p_{1}$ assigns a bundle with its first Pontrijagin class;

(3) via the Steenrod isomorphism, the homomorphism $i_{\ast }$
induced by the fibre inclusion corresponds to the operation $\zeta
\rightarrow \zeta \oplus \varepsilon ^{1}$, where $\oplus $ means
Whitney sum and where $\varepsilon ^{1}$ is the trivial $1$-bundle
over $S^{4}$; and

(4) the map triangle commutes by the stability of Pontrijagin classes.

From the Bott-periodicity we have $\pi _{3}(SO(5))=\mathbb{Z}$. It
is also known that $p_{1}:$ $Vect^{5}(S^{4})\rightarrow
H^{4}(S^{4};\mathbb{Z})$ is surjective onto the subgroup
$2H^{4}(S^{4};\mathbb{Z})\subset H^{4}(S^{4};\mathbb{Z})$.
Summarizing we have

(5) $p_{1}:$ $Vect^{5}(S^{4})\rightarrow H^{4}(S^{4};\mathbb{Z})$ is
injective.

From (5.2) we find that $p_{1}(\widehat{\xi })=0$. As a result (4) and (5)
imply that $i_{\ast }(\widehat{\xi })=0$. From the exactness of the top
sequence one concludes

$\qquad \widehat{\xi }=k\partial (\iota _{4})$ for some $k\in \mathbb{Z}$,

\noindent where $\iota _{4}\in \pi _{4}(S^{4})=\mathbb{Z}$ is the class of
identity map. Since $\partial (\iota _{4})=\tau $ (cf.[Wh, p.196]) we have $%
f(\widehat{\xi })=(2k,0)$ by the Example. This contradiction to (5-2)
completes the proof.$\square $

\textbf{Remark 4.} The proof of Lemma 9 indicates that the bundles $\tau $
and $\gamma $ in the Example generate the group $Vect^{4}(S^{4})$. These two
bundles were used by Milnor in [MS, p.247] to illustrate his original
construction of different differential structures on the $7$-sphere in 1956
[M].$\square $

\begin{center}
\textbf{References}
\end{center}

[B] A. Borel, Sur la cohomologie des espaces fibr\'{e}s principaux et des
espaces homog\`{e}nes de groupes de Lie compacts, Ann. of Math. (2)57(1953),
115-207.

[BS] R. Bott and H. Samelson, Application of the theory of Morse to
symmetric spaces, Amer. J. Math., Vol. LXXX, no. 4 (1958), 964-1029.

[D] H. Duan, Self-maps of the Grassmannian of complex structures, Compositio
Math. 132(2002), 159-175.

[HPT] W. Y. Hsiang, R. Palais and C. L. Terng, The topology of isoparametric
submanifolds, J. Diff. Geom., Vol. 27 (1988), 423-460.

[Hu] J. E. Humphreys, Introduction to Lie algebras and representation
theory, Graduated text in Math. 9, Springer-Verlag New York, 1972.

[LD] B. Li and H. Duan, Spin characteristic classes and reduced Kspin Groups
of low dimensional complex, Proc. Amer. Math. Soc., Vol. 113, No. 2 (1991),
479-491.

[M] J. Milnor, On manifolds homeomorphic to the 7-sphere. Ann. of Math. (2)
64 (1956), 399--405.

[MS] J. Milnor and J. D. Stasheff, Characteristic classes, Ann. Math.
Studies 76, Princeton University Press, Princeton NJ, 1974.

[S] N. A. Steenrod, The topology of fiber bundles, Princeton University
Press, 1951.

[Wh] G. W. Whitehead, Elements of Homotopy theory, Graduate texts in Math.
61, Springer-Verlag, New York Heidelberg Berlin, 1978.

\end{document}